\documentstyle[amsthm,amsfonts,amssymb,amsgen,amsmath,amsopn,verbatim]{article}

\theoremstyle{plain}
\newtheorem{thm}{Theorem}[section]
\newtheorem{lem}[thm]{Lemma}
\newtheorem{prop}[thm]{Proposition}

\theoremstyle{definition}
\newtheorem{rem}[thm]{Remark}
\newtheorem{eg}[thm]{Example}

\newcommand{\tT}{{\widetilde{T}}}
\newcommand{\wg}{\wedge}

\newcommand{\zero}{{\bf 0}}
\newcommand{\one}{{\bf 1}}
\newcommand{\bfu}{{\bf u}}
\newcommand{\bfv}{{\bf v}}
\newcommand{\obfi}{{\overrightarrow{\boldsymbol \i}}}
\newcommand{\obfj}{{\overrightarrow{\boldsymbol \j}}}
\newcommand{\obfk}{{\overrightarrow{\bf k}}}
\newcommand{\veps}{{\varepsilon}}
\newcommand{\bveps}{{\boldsymbol \varepsilon}}
\newcommand{\bfi}{{\bf i}}
\newcommand{\bfj}{{\bf j}}
\newcommand{\bfx}{{\bf x}}
\newcommand{\bfy}{{\bf y}}

\newcommand{\bfk}{{\bf k}}
\newcommand{\bfeps}{{\boldsymbol \varepsilon}}
\newcommand{\R}{{\mathbb R}}
\newcommand{\CP}{{\mathbb CP}}
\newcommand{\C}{{\mathbb C}}

\DeclareMathOperator{\id}{{id}}

\newcommand{\ga}{{\alpha}}
\newcommand{\gl}{{\lambda}}
\newcommand{\gb}{{\beta}}
\newcommand{\gd}{{\delta}}
\newcommand{\gep}{{\varepsilon}}
\newcommand{\eps}{{\epsilon}}
\newcommand{\gs}{{\sigma}}
\newcommand{\gTH}{{\Theta}}
\newcommand{\caL}{{\mathcal L}}
\newcommand{\caD}{{\mathcal D}}
\newcommand{\gemL}{{\mathfrak L}}
\newcommand{\SS}{{\mathfrak S}}
 
\begin{document}
 
\title{Analytic continuation of multiple polylogarithms
\footnote{Mathematics Subject Classification (2000): 
32D15, 53C65, 14D05, 33B30.}}
\author{Jianqiang Zhao\footnote{Partially supported by NSF grant DMS0139813}}
\date{}
\maketitle

\medskip
\noindent{\bf Abstract.}
In this paper we shall define the analytic
continuation of the multiple polylogarithms by using 
Chen's theory of iterated path integrals
and compute the monodromy of all multiple logarithms
explicitly.

\section{Introduction}
In recent years, there is a revival of interest in multi-valued
classical polylogarithms and their single-valued cousins. 
For any positive integer $m_1,\dots,m_n$,
Goncharov \cite{Gicm} generalizes the classical polylogarithms 
and defines the multiple  polylogarithms as follows:
\begin{equation}\label{pexp}
Li_{m_1,\dots,m_n}(x_1,\dots,x_n) =\sum_{0 <
k_1<k_2<\dots<k_n} \frac{x_1^{k_1}x_2^{k_2} \dots
x_n^{k_n}}{k_1^{m_1}k_2^{m_2}\dots k_n^{m_n}},\qquad |x_i|<1.
\end{equation}
Conventionally one refers $n$ as the {\em depth} and
$K:=m_1+\dots+m_n$ as the {\em weight}. When the depth
$n=1$ the function is nothing but the classical polylogarithm.
More than a century ago it was already known to H.~Poincar\'e \cite{Po}
that hyperlogarithms
$$F_n\left({{a_1,\dots, a_n}\atop{b_1,\dots,b_n}}\Big|z\right)
=\int_{b_n}^z\cdots\int_{b_2}^{t_3}\int_{b_1}^{t_2} \frac{d t_1}{t_1-
a_1}
\frac{d t_2}{t_2-a_2}\cdots \frac{d t_n}{t_n-a_n}$$
are important for solving differential equations. 
Notice that the multiple polylogarithm
\begin{equation}\label{old}
Li_{m_1,\dots,m_n}(x_1,\dots,x_n)=(-1)^n
F_K\bigg({{a_1,\overbrace{0,\dots,0}^{m_1-1\text{ times}},
\dots, a_n,\overbrace{0,\dots,0}^{m_n-1\text{ times}}}\atop
{0\ ,\ 0,\dots,0\ ,\dots,\ 0\ \, ,\ 0,\dots,0}}\Big|1\bigg),
\end{equation}
where $a_i=1/(x_i\dots x_n)$ for $1\le i\le n$. It is an
iterated path integral in the sense of Chen \cite{Ch1} whose
path lies in $\C$. One thus can 
easily enlarge its  domain of definition to some open subset of $\C^n$.
However, it is not obvious that this actually gives a genuine analytic
continuation in the usual sense. In particular, it is difficult to
study the monodromy of the multiple polylogarithms by this expression.
It is our primary goal in this paper to define the analytic 
continuation of the above function to $\C^n$ as a 
multi-valued meromorphic function by using Chen's iterated 
path integrals with all paths lying in $\C^n$.  

In early 1980s Deligne \cite{D} discovers that the dilogarithm gives 
rise to a good variation of mixed Hodge-Tate structures. 
This has been generalized to polylogarithms (cf.~\cite{Hain}) 
following Ramakrishnan's computation of the monodromy
of the polylogarithms. The monodromy computation also 
yields the single-valued variant $\caL_n(z)$ of the 
polylogarithms (cf.~\cite{BD,Zag}). These functions in turn
have significant applications in arithmetic such as Zagier's conjecture
\cite[p.622]{Zag}. On the other hand,
as pointed out in \cite{G,Gmzeta}, ``higher cyclotomy theory'' should 
study the multiple polylogarithm motives at roots of unity, not only
those of the polylogarithms. This is the primary reason for us
to focus our attention on multiple polylogarithms.

According to the theory of framed mixed Hodge-Tate structures
the multiple polylogarithms are period functions
of some variations of mixed Hodge-Tate structures (see \cite{BGSV},
\cite[\S12]{Gicm} and \cite[\S3.5]{Gicm}). 
However, in order to have ``reasonable'' variations 
we should be able to control their behavior at ``infinity'' $X_n$
(see \eqref{Xn}). Our analytic continuation 
of the multiple polylogarithms is more straight-forward
than \eqref{old} so that we can not only see that
the multiple polylogarithms are multi-valued but also  
determine $X_n$ explicitly where the singularities lie. 
We then compute the monodromy of all multiple logarithms 
$$\gemL_n(x_1,\dots,x_n):=
Li_{\underbrace{\scriptstyle 1,\dots,1}_{\scriptstyle n\text{ times}}}
(x_1,\dots,x_n).$$
This enables us to calculate mixed Hodge structures associated with
some families of multiple polylogarithms in another paper \cite{Zmhs}
including all the multiple logarithms.

We point out that the mixed Hodge structures of iterated integrals over
$\CP\setminus\{0,1,\infty\}$ have been studied by Wojtkowiak \cite{Wo}.
His approach is completely different from ours.

The motivation of this paper comes from \cite[\S2,3]{G} 
where the Hodge-Tate structures associated with the
double logarithms are discussed. The author thanks
his advisor Sasha Goncharov for his constant encouragement
and Herbert Gangl for kindly informing the author of
the preprint \cite{Wo} of Wojtkowiak.
 
\section{Preliminaries on Chen's theory of iterated path integrals} 
The main references of this section
are two of Chen's papers \cite{Ch1} and \cite{Ch2}.

For a 1-form $f(t)dt$ over $\R$ the integral $\int_a^b f(t)dt$
is understood in the usual way. For $r>1$, define inductively
$$\int_a^b f_1(t)dt\cdots  f_r(t)dt=\int_a^b
\left(\int_a^t f_1(\tau )d\tau \cdots  f_{r-1}(\tau )d\tau \right)
f_r(t) dt.$$
When $r=0$, set the integral to be 1. For example, the classical
polylogarithm
$$Li_n(x)= \int_0^x \frac{dt}{1-t}
\underbrace{\frac{dt}{t}\cdots\frac{dt}{t}}_{(n-1)\text{ times}}.$$

More generally, let $w_1,w_2,\dots$ be 1-forms on a manifold $M$ 
and let $\ga:[0,1]\rightarrow M$ be a piecewise smooth path. Write
$$\ga^*w_i=f_i(t)dt$$
and define the {\em iterated path integral}
$$\int_\ga w_1\cdots w_r=\int_0^1 f_1(t)dt\cdots  f_r(t)dt.$$

The following results are crucial for the application of the Chen's
theory of iterated path integrals.
\begin{lem}\label{chen's}
Let $w_i$ $(i\ge 1)$ be $\C$-valued 1-forms on a manifold $M$.

(i) The value of $\int_\ga w_1\cdots w_r$ is independent of the
parameterization of $\ga$.

(ii) If $\ga,\gb:[0,1]\longrightarrow M$ are composable paths
(i.e. $\ga(1)=\gb(0)$), then
$$\int_{\ga\gb} w_1\cdots w_r=\sum_{j=0}^r
\int_\ga w_1\cdots w_i\int_\gb w_{i+1}\cdots w_r.$$
Here, we set $\int_\ga \phi_1\cdots \phi_m=1$ if $m=0$.

(ii) For every path $\ga$,
$$ \int_{\ga^{-1}}w_1\cdots w_r = (-1)^r\int_\ga w_r\cdots w_1.$$

(iv) For every path $\ga$,
$$ \int_\ga w_1\cdots w_r\int_\ga w_{r+1}\cdots w_{r+s}
=\sum_\gs \int_\ga w_{\gs(1)}\cdots w_{\gs(r+s)}$$
where $\gs$ ranges over all shuffles of type $(r,s)$, i.e., 
permutations $\gs$ of $r+s$ letters with 
$\gs^{-1}(1)<\cdots<\gs^{-1}(r)$
and $\gs^{-1}(r+1)<\cdots< \gs^{-1}(r+s)$.
\end{lem}
\begin{proof}
(i) can be derived from
the theorem on \cite[p.~361]{Ch1}. (ii) and (iii) are formulas
(1.6.1) and (1.6.2) of \cite{Ch1} respectively. Ree \cite{Ree}
discovered the shuffle relation (iv) which appeared as (1.5.1) 
in \cite{Ch1}.
\end{proof}
\begin{lem}\label{pathindpt}
If $w_i^{(j)}$ are closed 1-forms for $1\le i\le r$ and 
$1\le j\le n$ such that 
$\sum_{j} w_1^{(j)}\wg w_2^{(j)}=\sum_{j} w_2^{(j)}\wg w_3^{(j)}
=\cdots=\sum_{j} w_{r-1}^{(j)}\wg w_r^{(j)}=0$ then
$\sum_{j} \int_\ga w_1^{(j)} w_2^{(j)}\cdots w_r^{(j)}$ 
only depends on the homotopy class of $\ga$.
\end{lem}
\begin{proof} The case $j=1$ is proved on \cite[p.~366]{Ch1}. The
case $r=2$ can be found on \cite[p.~368]{Ch1}. The general case
follows from a similar argument.
\end{proof}

\section{The index set $\SS(m_1,\dots,m_n)$}
We first introduce an index set with two 
different kinds of orderings. 
 
\medskip
\noindent{\bf 2.1. Definition.} 
Define the index set
$$\SS(m_1,\dots,m_n)=\{\bfi=(i_1,\dots,i_n):\ 0\le i_t\le m_t
\text{ for }t=1,\cdots, n\}$$
and the weight function $| \cdot |$  on $\SS(m_1,\dots,m_n)$ by
$$|(i_1,\dots,i_n)|=i_1+\cdots+i_n.$$
For brevity, we write 
$\zero=(0,\dots,0)\in \SS(m_1,\dots,m_n)$
which is the only index of weight $0$ in $\SS(m_1,\dots,m_n)$
and  $\one_K=(m_1,\dots,m_n)\in \SS(m_1,\dots,m_n)$ which is the
only index of the highest weight $K:=m_1+\cdots+m_n$ in
$\SS(m_1,\dots,m_n)$. We also define the depth function of 
the index $(i_1,\dots,i_n)$ by $\sharp\{t:\, i_t\ne 0\},$ 
i.e., the number of nonzero components.

We shall need two orderings on $\SS(m_1,\dots,m_n)$.

\medskip
\noindent{\bf 2.2. A complete ordering.} 
The complete ordering is defined as follows. Let 
$\bfi=(i_1,\dots,i_n)$ and $\bfj=(j_1,\dots,j_n).$
If $|\bfi|<|\bfj|$ then $\bfi<\bfj$ (or, equivalently, $\bfj>\bfi$).
If $|\bfi|=|\bfj|$ then $\bfi>\bfj$ if $\max\{i_t: 1\le t\le n\}>
\max\{j_t: 1\le t\le n\}$. Otherwise, we compare the second largest 
components
of $\bfi$ and $\bfj$, and so on. If
$\{i_1,\dots,\i_n\}=\{j_1,\dots,j_n\}$ as two set
then the usual lexicographic order from left to right
is in force with $0<1<\cdots $.
For instance, $(0,0,1)<(1,0,1)<(1,1,0)<(0,2,0)$ in $\SS(1,2,1)$.
\begin{rem} In the multiple logarithm case, namely, when
$m_1=\dots=m_n=1$, there is a one-to-one correspondence between
$\SS_n$ and the set of non-negative integers less than
$2^n$. Thus one is tempted to use the
conventional order of positive integers in binary forms.
However this is not suitable in our situation.
\end{rem}

\medskip
\noindent{\bf 2.3. A partial ordering and the retraction map.} 
The partial ordering is defined as follows.
Let $\bfi=(i_1,\dots,i_n)$
and $\bfj=(j_1,\dots,j_n)$. We set $\bfj\prec \bfi$
(or, equivalently, $\bfi \succ\bfj$) if $j_t\le i_t$ for every 
$1\le t\le n$. For example $(0,0,1,0)\prec (0,1,1,0)$ in 
$\SS(1,1,1,1)$ but $(1,0,0,0)\not\prec (0,1,1,0)$ and 
$(1,0,0,0)\not\succ (0,1,1,0)$. Clearly
$\bfj\prec \bfi$ implies $\bfj< \bfi$ but not vice versa.

Suppose $\bfi$ has depth $k$ with
$i_{\tau_s}\ne 0$ for $1\le s\le k$ while $\bfj$ has depth $l$ and
$j_{t_r}\ne 0$ for $1\le r\le l$. If $\bfj\prec \bfi$ then 
we can write $t_r=\tau_{\ga_r}$ for $1\le r\le l$.
For such $\bfi$ and $\bfj$ we define the $\bfi$-th
retraction map $\rho_\bfi$ from $\SS(m_1,\dots,m_n)$ to
$\SS(i_{\tau_1},\dots,i_{\tau_k})$ as follows. The entry
of $\rho_\bfi(\bfj)$ is $j_{\tau_{\ga_r}}$ if it is
at the $\ga_r$-th ($1\le r\le l$) component and 0 at all 
other components.
For instance $\rho_{(02010)}\big((01000)\big)=(10)\in \SS(2,1)$.
In particular,  $\rho_\bfi(\bfi)=(i_{\tau_1},\dots,i_{\tau_k})$
has highest weight in $\SS(i_{\tau_1},\dots,i_{\tau_k})$.

\medskip
\noindent{\bf 2.4. Vector indices.} 
Let $\SS^K(m_1,\dots,m_n)$ be the set of $K$-tuples
$\obfj=(\bfj_1,\dots,\bfj_K)$ of $\SS(m_1,\dots,m_n)$
such that $|\bfj_t|=t$
and  $\bfj_1\prec \cdots \prec\bfj_K=\one_K$. One may think $\obfj$
as a length $K$ queue of indices of $\SS(m_1,\dots,m_n)$ in which
each index is produced by increasing some component of the preceding
index by 1.

\medskip
\noindent{\bf 2.5. Additional notation.} 
Throughout the paper we fix 
$\bfu_s:=(0,\dots,0,1,0,\dots,0)\in \SS(m_1,\dots,m_n)$ of
weight $1$ where the entry $1$ is at the $s$-th component. Whenever
the  $s$-th component $i_s$ of $\bfi$ satisfies $i_s<m_s$ we can
increase $i_s$ by 1 to get a new index which is denoted by
$\bfi+\bfu_s$. If $i_s>0$ we similarly define $\bfi-\bfu_s$
as the index with the $s$-th component of $\bfi$ decreased by 1.
Fix $\bfv_s=\one_K-m_s\bfu_s\in \SS(m_1,\dots,m_n)$ 
whose components are nonzero except at the $s$-th position.

When $m_1=\cdots=m_n=1$ we write $\SS(1,\dots,1)=\SS_n$.

\medskip
\noindent{\bf 2.6. Transposition functions.}
Fix arbitrary 
$\obfj=(\bfj_1,\dots,\bfj_K)\in \SS^K(m_1,\dots,m_n)$ 
and $1<r\le K$ we write
$$\bfj_r=\bfj_{r-1}+\bfu_s=(t_1,\dots,t_s,0,\dots,0,t_a,\dots,t_n),
\quad 0\le s<a\le n+1,\ t_a\ne 0.$$
Here if $a=n+1$ then the last nonzero component of $\bfj_r$ is $t_s$.
We define the transposition functions on 
$\bfi=(i_1,\dots,i_n)\in\SS(m,\dots,m)$
with $m=\max\{m_1,\dots,m_n\}$ by
$$T^r_0=\id, \quad T^r_1 (\bfi)=(i_{\gs(1)},\dots,i_{\gs(n)})$$
where if $t_s>1$ or $a=n+1$ then $\gs=\id$ whereas if $t_s=1$ 
and $a\le n$ then $\gs$ is the
transposition in the symmetric group of $n$ elements that exchanges
$s$ and $a$.

\section{Analytic continuation of multiple polylogarithms}
Let $\bfx=(x_1,\dots,x_n)$ be a variable over $\C^n$. Define
$$S_n=\C^n\setminus X_n,\quad S'_n=\C^n\setminus X'_n,$$
where the divisors are defined by
\begin{align}\label{Xn}
X_n=&\Bigl\{\bfx\in\C^n: \prod_{1\le i\le n}(1-x_j)
\prod_{1\le j<k\le n} \big(1-x_j\dots x_k\big)=0\Bigr\},\\
X'_n=&\Bigl\{\bfx\in\C^n: \prod_{1\le i\le n}x_j(1-x_j)
\prod_{1\le j<k\le n} \big(1-x_j\dots x_k\big)=0\Bigr\}. \notag
\end{align}
It is clear that $S'_n\subset S_n$. Set
$$D_n=\Bigl\{(x_1,\dots,x_n)\in \C^n: 
	\Bigl|x_j-\frac 12\Bigr|<\frac 12,
	j=1,\dots, n\Bigr\}\subset S'_n.$$
Denote a varying base point by 
$\bveps_n=(\veps,\dots,\veps)\in D_n$.

Suppose the depth of $\bfi=(i_1,\dots,i_n)$ is $k$ and
$i_{\tau_1}\ne 0,\dots,i_{\tau_k}\ne 0$. We define
$$a_t=a_t(\bfx):=(x_t\dots x_n)^{-1} \text{ for }1\le t\le n$$
and
\begin{equation}\label{y's}
\bfx(\bfi)=\bfy=(y_1,\dots,y_k),\qquad
y_m=\prod_{\ga=\tau_m}^{\tau_{m+1}-1}x_\ga
	=\frac{a_{\tau_{m+1}}(\bfx)}{a_{\tau_m}(\bfx)}, 
	\quad 1\le m\le k
\end{equation}
with $\tau_{k+1}=n+1$ and $a_{n+1}=1$. We also write
$a_m(\bfy)=(y_m\dots y_k)^{-1}=a_{\tau_m}(\bfx)$. Note that 
$\bfx(\bfi)\in \C^k$ which is the reason 
why we call $k$ the depth of $\bfi$.

We begin with some 1-forms which will be used to express
the multiple polylogarithms. Take
$\obfj\in \SS^K(m_1,\dots,m_n)$ and $\bfj_r=\bfj_{r-1}+\bfu_s$
as given in \S2.6. 
For any $(\gd_1,\dots,\gd_K)\in \SS_K$, namely, $\gd_t=0$ or 1,
let $\bfy=\bfx$ if $r=K$ and 
$$\bfy=(y_1,\dots,y_l)=\bfx 
\big(T_{\gd_{r+1}}^{r+1}\circ \cdots\circ T_{\gd_K}^K(\bfj_r) \big) 
\qquad \text{ if }\quad 1\le r<K$$ 
where $l$ is the depth of $\bfj_r$ because the transposition 
functions do not change the depth of an index.
We let $t_{\ga_1}\ne 0,\dots,t_{\ga_l}\ne 0$ and $s=\ga_\gl$
(because $t_s\ne 0$) and set 
$$w_{\obfj}^{r,\gd_r}(\bfy):=\begin{cases}
0 \quad &\text{ if $t_s>1$ and $\gd_r=1$},\\
dy_\gl/ y_\gl \quad &\text{ if $t_s>1$ and $\gd_r=0$},\\
dy_\gl/(1-y_\gl) \quad &\text{ if $t_s=1$ and $\gd_r=0$} ,\\
dy_\gl/y_\gl(y_\gl-1) \quad
	&\text{ if $\gl<l$, $t_s=1$ and $\gd_r=1$},\\
0 \quad &\text{ if $\gl=l$, $t_s=1$ and $\gd_r=1$}.
\end{cases}$$
It is obvious that $w_{\obfj}^{r,\gd_r}(\bfy)$ is always 
a closed 1-form whose singularities lie only along $X'_n$.

\begin{prop}\label{genL} Let $\int_p \bigsqcup_{r=1}^K w_r$ 
denote the iterated integral $\int_p w_1\cdots w_K$. Then
for every $\bfx\in D_n$
$$Li_{m_1,\dots,m_n}(\bfx)=\lim_{\veps\to 0} \int_{\bveps_n}^{\bfx}
\sum_{\substack{(\gd_1,\dots,\gd_K)\in \SS_K\\
\obfj\in\SS^K(m_1,\dots,m_n)}} \operatorname*{\bigsqcup}_{r=1}^K
w_{\obfj}^{r,\gd_r}\left(\bfx\Bigl(T_{\gd_{r+1}}^{r+1}
	\circ \cdots\circ T_{\gd_K}^K(\bfj_r) \Bigr)\right)$$
where the paths of the iterated integral lie entirely in $D_n$
\end{prop}
\begin{proof} 
We prove this by induction on $K$. When $K=1$ this is
trivial.
Assume $K>1$ and the proposition is true for $K-1$.
Using the power series expansion \eqref{pexp} 
it is straight-forward to check that
$$dLi_{m_1,\dots,m_n}(\bfx)
	=\sum_{t=1}^n d_tLi_{m_1,\dots,m_n}(\bfx)$$
where if $m_t>1$ then
\begin{equation}\label{dm2}
d_tLi_{m_1,\dots,m_n}(\bfx)
=Li_{m_1,\dots,m_{t-1},m_t-1,m_{t+1},\dots,m_n}(\bfx) dx_t/ x_t
\end{equation}
whereas if $m_t=1$ then
\begin{multline}\label{dtL}
d_tLi_{m_1,\dots,m_n}(\bfx)=
Li_{m_1,\dots,m_{t-1},m_{t+1},\dots,m_n}\big(\bfx(\bfv_t)\big)
 dx_t/(1-x_t)\\
+Li_{m_1,\dots,m_{t-1},m_{t+1},\dots,m_n}\big(\bfx(\bfv_{t+1})\big)
dx_t/x_t(x_t-1).
\end{multline}
Here when $t=n$ and $m_n=1$ the second term in the sum 
does not appear. Observe that for any 
$\obfj\in \SS^K(m_1,\dots,m_{t-1},1,m_{t+1},\dots,m_n)$ with
$\bfj_{K-1}= (m_1,\dots,m_{t-1},0,m_{t+1},\dots,m_n)$ and $t<n$ 
we have  $T_1^{K}(\bfj_{K-1})=
(m_1,\dots,m_{t-1},m_{t+1},0,m_{t+2},\dots,m_n)$ 
and therefore
$$\bfx(\bfv_{t+1})=\bfx(m_1,\dots,m_t,0,m_{t+2},\dots,m_n)
=\bfx\big(T_1^{K}(\bfj_{K-1})\big).$$
Hence
\begin{equation}\label{ddd}
dLi_{m_1,\dots,m_n}(\bfx)= 
\sum_{\gd_K=0,1}
\sum_{\substack {|\bfj_{K-1}|=K-1,\\ 
	\bfj_{K-1}\in\SS(m_1,\dots,m_n) }} 
Li_{T_{\gd_K}^K(\bfj_{K-1})}
\Bigl(\bfx\big(T_{\gd_K}^K(\bfj_{K-1})\big)\Bigr)
w_{\obfj}^{K,\gd_K}(\bfx)
\end{equation}
where we write $Li_{\bfv_t}=Li_{m_1,\dots,m_{t-1},m_{t+1},\dots,m_n}.$

For $1\le t\le n$ define the embeddings in the obvious way
$$
\iota_s:\SS(m_1,\dots,m_{s-1},m_s-1,m_{s+1},\dots,m_n)
	\hookrightarrow\SS(m_1,\dots,m_n)
$$
where when $m_s=1$ the left hand side is understood as
$\SS(m_1,\dots,m_{s-1},m_{s+1},\dots,m_n)$ which is 
identified with $\SS(m_1,\dots,m_{s-1},0,m_{s+1},\dots,m_n)$ 
as a subset of the right hand side.
By abuse of notation, we further define
$$\aligned
\iota_s:\SS^{K-1}(m_1,\dots,m_{s-1},m_s-1,m_{s+1},\dots,m_n)&
\longrightarrow \SS^K(m_1,\dots,m_n)\\
(\bfk_1,\dots,\bfk_{K-1})&\longmapsto (\bfk_1,\dots,\bfk_{K-1},\one_K).
\endaligned$$

By induction 
\begin{multline*}
Li_{T_{\gd_K}^K(\bfj_{K-1})}
\Bigl(\bfx\big(T_{\gd_K}^K(\bfj_{K-1})\big)\Bigr)
=\lim_{\veps\to 0}\int_{\bveps}^
{\bfx{\textstyle (}T_{\gd_K}^K(\bfj_{K-1}){\textstyle )}}
\sum_{(\gd_1,\dots,\gd_{K-1})\in \SS_{K-1}}
\sum_{\obfk\in\SS^{K-1}{\textstyle (}
	T_{\gd_K}^K(\bfj_{K-1}){\textstyle )}} \\
\operatorname*{\bigsqcup}_{r=1}^{K-1}
w_{\obfk}^{r,\gd_r}\left(\bfx\big(T_{\gd_K}^K(\bfj_{K-1})\big)
\Bigl(\tT_{\gd_{r+1}}^{r+1}
	\circ \cdots\circ \tT_{\gd_{K-1}}^{K-1}(\bfk_r) \Bigr)\right)
\end{multline*}
where the transposition functions $\tT$ may differ from $T$.
Let $\bfj_K=\bfj_{K-1}+\bfu_s$. We now show that
\begin{equation}\label{keycomm}
\bfx\big(T_{\gd_K}^K(\bfj_{K-1}) \big)\Bigl(\tT_{\gd_{r+1}}^{r+1}
	\circ \cdots\circ \tT_{\gd_{K-1}}^{K-1}(\bfk_r) \Bigr)
=\bfx\Bigl(T_{\gd_{r+1}}^{r+1}
	\circ \cdots\circ T_{\gd_K}^K\big(\iota_s(\bfk_r)\big) \Bigr).
\end{equation}
This is trivial if $\gd_K=0$ or $s=n$ or $t_s>1$ where $t_s$ is the
$s$-th component of $\bfj_K$ because in these cases we have
$T_{\gd_K}^K=\id$, $\tT=T$ and $\iota_s=\id$. So we assume $\gd_K=1$, 
$t_s=1$ and $s<n$. Then
$$\bfx\big(T_{\gd_K}^K(\bfj_{K-1}) \big)=
(x_1,\dots,x_{s-1},x_sx_{s+1},x_{s+2},\dots,x_n),$$
and the $s$-component of $\iota_s(\bfk_r)$ is $0$ 
by definition. By straight-forward computation we find that
$$T_{\gd_K}^K\circ \iota_s\circ \tT_{\gd_\gl}^{\gl}
=T_{\gd_\gl}^\gl\circ T_{\gd_K}^K\circ 
\iota_s\qquad \text{ for all } r+1\le \gl\le K-1.$$
These implies equation \eqref{keycomm} immediately because
$$\aligned
\bfx\big(T_{\gd_K}^K(\bfj_{K-1}) \big)\Bigl(\tT_{\gd_{r+1}}^{r+1}
	\circ \cdots\circ \tT_{\gd_{K-1}}^{K-1}(\bfk_r) \Bigr)
=&\bfx\Bigl(T_{\gd_K}^K\circ\iota_s\circ\tT_{\gd_{r+1}}^{r+1}
	\circ \cdots\circ \tT_{\gd_{K-1}}^{K-1}(\bfk_r) \Bigr)\\
=&\bfx\Bigl(T_{\gd_{r+1}}^{r+1}
	\circ \cdots\circ T_{\gd_K}^K\big(\iota_s(\bfk_r)\big) \Bigr).
\endaligned$$
Therefore from \eqref{ddd} and the one-to-one correspondence:
$$\iota_s:\SS^{K-1}\big(T_{\gd_K}^K(\bfj_{K-1}) \big)
\longleftrightarrow \big\{(\bfi_1,\dots,\bfi_K)\in 
\SS^K(m_1,\dots,m_n): \bfi_{K-1}=T_{\gd_K}^K(\bfj_{K-1})\big\}$$
we see that $dLi_{m_1,\dots,m_n}(\bfx)$ is equal to
$$
\sum_{\substack{(\gd_1,\dots,\gd_K)\in \SS_K\\
	\obfj\in\SS^K(m_1,\dots,m_n)}} 
\left(\lim_{\veps\to 0}\int_{\bveps}^
{\bfx{\textstyle (}T_{\gd_K}^K(\bfj_{K-1}){\textstyle )}}
\operatorname*{\bigsqcup}_{r=1}^{K-1}
w_{\obfj}^{r,\gd_r}\left(\bfx\Bigl(T_{\gd_{r+1}}^{r+1}
	\circ \cdots\circ T_{\gd_K}^K(\bfj_r) \Bigr)\right)\right)
	w_{\obfj}^{K,\gd_K}(\bfx).
$$
This finishes the proof of the proposition by induction because
$\lim_{\veps\to 0} Li_{m_1,\dots,m_n}(\bveps_n)=0$
where the limiting process takes place inside $D_n$.
\end{proof}

By the above proposition we can define the analytic 
continuation of $Li_{m_1,\dots,m_n}(\bfx)$ to $S'_n$ as the 
iterated path integral
\begin{equation}\label{e:genL}
Li_{m_1,\dots,m_n}(\bfx)=\lim_{\veps\to 0} \int_{\bveps_n}^{\bfx}
\sum_{\substack{(\gd_1,\dots,\gd_K)\in \SS_K\\
\obfj\in\SS^K(m_1,\dots,m_n)}} \operatorname*{\bigsqcup}_{r=1}^K
w_{\obfj}^{r,\gd_r}\left(\bfx\Bigl(T_{\gd_{r+1}}^{r+1}
	\circ \cdots\circ T_{\gd_K}^K(\bfj_r) \Bigr)\right),
\end{equation}
where all the paths lie inside $S'_n$.
Note that all the $1$-forms appearing in \eqref{e:genL} are
rational forms with logarithmic singularities along $X'_n$.

\begin{eg}\label{egss}
When $n=1$,
$$Li_1(x)=\int_0^x d\log\Bigl(\frac{1}{1-x}\Bigr)=-\log(1-x).$$
When $n=2$, $\SS_2=\{(0,0),(0,1),(1,0),(1,1)\}$ and there are two
elements in $\SS^2$: $((0,1),(1,1))$ and $((1,0),(1,1))$.
Let $\bfx=(x,y)$ then
$$\bfx(0,1)=y,\quad \bfx(1,0)=xy,\quad \bfx(1,1)=(x,y).$$
Thus
$$\aligned
Li_{1,1}(x,y)=&\int_{\zero}^{\bfx}
w_1(\bfx(0,1))w_1(\bfx)+ w_1(\bfx(1,0))w_2(\bfx)\\
=&\int_{(0,0)}^{(x,y)}
\frac{dy}{1-y} \frac{dx}{1-x}+\frac{d(xy)}{1-xy}
\left(\frac{dy}{1-y}+\frac{dx}{x(x-1)}\right).
\endaligned$$
When $n=3$ let $\bfx=(x,y,z)$. Then
$$\aligned
\bfx(0,0,1)=&z,\quad \bfx(0,1,0)=yz,\quad \bfx(1,0,0)=xyz,\quad
\bfx(0,1,1)=(y,z),\\
\bfx(1,0,1)=&(xy,z),\quad \bfx(1,1,0)=(x,yz),
	\quad \bfx(1,1,1)=(x,y,z).
\endaligned$$
Thus
$$\aligned
Li_{1,1,1}(x,y,z)=&\int_{(0,0,0)}^{(x,y,z)}
\frac{dz}{1-z}\frac{dy}{1-y} \frac{dx}{1-x}
+\frac{d(yz)}{1-yz}
\left(\frac{dz}{1-z}+\frac{dy}{y(y-1)}\right)
\frac{dx}{1-x}\\
+&\frac{d(yz)}{1-yz}\frac{dx}{1-x}
\left(\frac{dz}{1-z}+\frac{dy}{y(y-1)}\right)
+\frac{dz}{1-z}\frac{d(xy)}{1-xy}
\left(\frac{dy}{1-y}+\frac{dx}{x(x-1)}\right)\\
+&\frac{d(xyz)}{1-xyz}
\left(\frac{dz}{1-z}+\frac{d(xy)}{xy(xy-1)}\right)
\left(\frac{dy}{1-y}+\frac{dx}{x(x-1)}\right)\\
+& \frac{d(xyz)}{1-xyz}
\left(\frac{d(yz)}{1-yz}+\frac{dx}{x(x-1)}\right)
\left(\frac{dz}{1-z}+\frac{dy}{y(y-1)}\right).
\endaligned$$
\end{eg}

\begin{lem} \label{indpath}
The iterated path integral inside the limit of \eqref{e:genL} 
depends only on the homotopy class of the path 
from $\bveps_n$ to $\bfx$.
\end{lem}
\begin{proof} We use induction on the weight $K$ and
Lemma \ref{pathindpt} to prove this lemma.

When $K=1$ this is trivial. When $K=2$ there are two possibilities:
the dilogarithm $Li_2(x)$ and the double logarithm 
$Li_{1,1}(x_1,x_2)$. First
$$Li_2(z)=\int_0^z \frac{dx}{1-x}\frac{dx}{x}$$
and thus $\frac{dx}{1-x}\wg \frac{dx}{x}=0$ over $\C$.
Second, for the double logarithm $Li_{1,1}(x,y)$ as given
in Example~\ref{egss} we clearly have
\begin{equation}\label{keyrel}
\frac{dy}{1-y}\wg \frac{dx}{1-x}+
\frac{d(xy)}{1-xy}\wg
\left(\frac{dy}{1-y}-\frac{dx}{1-x}-\frac{dx}{x}\right)=0.
\end{equation}
Suppose now $K\ge 3$ and the lemma is proved for up to $K-1$.
Then it is not hard to see that we only need to show that the sum of
the wedge products of the last two 1-forms is zero. Let us 
look at the 2-form $dx_s\wg dx_t$ for $1\le s\ne t \le n$.
If $|s-t|>1$ then by the symmetry
of the equations \eqref{dm2} and \eqref{dtL} with respect 
to $s$ and $t$ and skewsymmetry of the wedge product 
we get the desired result. So we may assume that $|t-s|=1$.

\begin{description}
\item[{\sl{\em (i)}}] If $m_s>1$ and $m_t>1$ then without loss of 
generality
we may assume that $t=s+1$. We have
\begin{align}\label{ds1}
d_sLi_{m_1,\dots,m_n}(\bfx)=&
Li_{m_1,\dots,m_{s-1},m_s-1,m_{s+1},\dots,m_n}(\bfx)
	\frac{dx_s}{x_s}\\
d_{s+1}Li_{m_1,\dots,m_n}(\bfx)=&
Li_{m_1,\dots,m_s,m_{s+1}-1,m_{s+2},\dots,m_n}(\bfx)
\frac{dx_{s+1}}{x_{s+1}}. \label{dt1}
\end{align}
Hence, by skewsymmetry of the wedge product, the sum of 
$(d x_s\wg d x_{s+1})$-terms cancel
with the sum of $(d x_{s+1}\wg d x_s)$-terms.

\item[{\sl{\em (ii)}}] If $m_s=1$ and $m_t>1$ then we have \eqref{dt1}
and \eqref{dtL} with $t$ replaced by $s$ and get
$$\frac {dx_s}{x_s(x_s-1)}\wg \frac {dx_t}{x_t}
+\frac {d(x_sx_t)}{x_sx_t}\wg \frac {dx_s}{x_s(x_s-1)}
+\frac {dx_s}{1-x_s}\wg \frac {dx_t}{x_t}
	+\frac {dx_t}{x_t}\wg \frac{dx_s}{1-x_s}=0.
$$
Here if $s=n$ then the first two terms do not occur.

\item[{\sl{\em (iii)}}] If $m_s=m_t=1$ then we may assume that $t=s+1$.
Take \eqref{dtL} with $t$ replaced by $s$
and $s+1$ respectively. Quickly we find that the sum of
$(d x_s\wg d x_{s+1})$-terms is
\begin{multline*}
\frac {dx_{s+1}}{1-x_{s+1}}\wg \frac{dx_s}{1-x_s}
+\frac {dx_{s+1}}{x_{s+1}(x_{s+1}-1)}\wg \frac{dx_s}{1-x_s}\\
+\frac {d(x_sx_{s+1})}{x_sx_{s+1}(x_sx_{s+1}-1)}\wg 
	\frac {dx_s}{x_s(x_s-1)}
+\frac {d(x_sx_{s+1})}{1-x_sx_{s+1}}\wg \frac {dx_s}{x_s(x_s-1)}
=\frac {dx_{s+1}}{x_{s+1}}\wg  \frac{dx_s}{x_s}
\end{multline*}
whereas by symmetry the sum of $(d x_{s+1}\wg d x_s)$-terms is
$$ \frac{dx_s}{x_s}\wg \frac {dx_{s+1}}{x_{s+1}}$$
which cancels with the sum of the $(d x_s\wg d x_{s+1})$-terms.
\end{description}

By induction the lemma now follows from Lemma \ref{pathindpt}.
\end{proof}

We now show that the multiple polylogarithms have trivial monodromy
about each $x_j=0$, $j=1,\cdots,n,$ and therefore
they are actually well-defined on $S_n$. 
\begin{thm}\label{thm:mono}
Let $p(\veps)$ be a path in from $\bveps_n\in D_n$ to an arbitrary
$\bfx\in S'_n$.
Let $q(\veps)$ be a loop in $S'_n$ based at $\bveps_n\in D_n$
around any $\caD_{j0}=\{x_j=0\}$ ($j=1,\cdots,n$)
but no other irreducible components of $X'_n$, then
\begin{equation}\label{monoKcase}
\lim_{\veps\to 0}\left(\int_{q(\veps)p(\veps)}-\int_{p(\veps)}\right)
\sum_{\substack{(\gd_1,\dots,\gd_K)\in \SS_K\\
\obfj\in\SS^K(m_1,\dots,m_n)}} \operatorname*{\bigsqcup}_{r=1}^K
w_{\obfj}^{r,\gd_r}\left(\bfx\Bigl(T_{\gd_{r+1}}^{r+1}
	\circ \cdots\circ T_{\gd_K}^K(\bfj_r) \Bigr)\right)=0.
\end{equation}
Therefore the multiple polylogarithm $Li_{m_1,\dots,m_n}(\bfx)$
is a multi-valued holomorphic
function on $S_n$ and can be expressed by
\begin{equation}\label{defn0}
Li_{m_1,\dots,m_n}(\bfx)=\int_{\zero}^{\bfx}
\sum_{\substack{(\gd_1,\dots,\gd_K)\in \SS_K\\
\obfj\in\SS^K(m_1,\dots,m_n)}} \operatorname*{\bigsqcup}_{r=1}^K
w_{\obfj}^{r,\gd_r}\left(\bfx\Bigl(T_{\gd_{r+1}}^{r+1}
	\circ \cdots\circ T_{\gd_K}^K(\bfj_r) \Bigr)\right).
\end{equation}
\end{thm}
\begin{proof} We prove \eqref{monoKcase} in the lemma by induction on 
$K$. 
If $K=1$ clearly $dx_1/(1-x_1)$ has no singularity at 
$\caD_{10}=\{x_1=0\}$ and \eqref{defn0} is
obvious. Assume the cases up to $K-1$ are true. Consequently, 
if $m_1+\dots+m_n=K-1$ then $Li_{m_1,\dots,m_n}(\bfx)$ 
is well defined by \eqref{defn0} and is equal to 0 if any $x_i=0$
and the path from $\zero$ to $\bfx$ does not enclose any irreducible
component of $X_n$.

We first prove that
\begin{equation}\label{Kcase}
\lim_{\veps\to 0} \int_{q(\veps)} 
\sum_{\substack{(\gd_1,\dots,\gd_K)\in \SS_K\\
\obfj\in\SS^K(m_1,\dots,m_n)}} \operatorname*{\bigsqcup}_{r=1}^K
w_{\obfj}^{r,\gd_r}\left(\bfx\Bigl(T_{\gd_{r+1}}^{r+1}
	\circ \cdots\circ T_{\gd_K}^K(\bfj_r) \Bigr)\right)=0.
\end{equation}
This is a special case of \eqref{monoKcase} when $p(\veps)$ 
shrinks to a point.

By equation \eqref{ddd} in the proof in Proposition \ref{genL}
we find that the left hand side of \eqref{Kcase} is equal to
$$\lim_{\veps\to 0} \int_{q(\veps)}  \sum_{\gd_K=0,1}
\sum_{\substack {|\bfj_{K-1}|=K-1,\\ 
	\bfj_{K-1}\in\SS(m_1,\dots,m_n) }} 
Li_{T_{\gd_K}^K(\bfj_{K-1})}
\Bigl(\bfx\big(T_{\gd_K}^K(\bfj_{K-1})\big)\Bigr)
w_{\obfj}^{K,\gd_K}(\bfx).$$
By Lemma \ref{indpath} we may assume that the path 
$q(\veps)$ lies in the (real) two dimensional space 
$\cap_{i\ne j}\{x_i=\veps\}$
enclosing $\caD_{j0}$ clockwise only once. By induction 
the terms with $\bfj_K-\bfj_{K-1}\ne \bfu_j$ in the above sum 
are clearly zero. Thus the integral is reduced to
\begin{equation}\label{iszero}
\lim_{\veps\to 0} \int_{q(\veps)}  \sum_{\gd_K=0,1}
  Li_{T_{\gd_K}^K(\bfj_{K-1})}
	\Bigl(\bfx\big(T_{\gd_K}^K(\bfj_{K-1})\big)\Bigr)
	w_{\obfj}^{K,\gd_K}(\bfx)),\qquad \bfj_{K-1}=\one_K-\bfu_j
\end{equation}
since $\bfj_K=\one_K$. By the induction assumption the 
functions in front of $w_{\obfj}^{K,\gd_K}(\bfx)$ are 
regular along $\caD_{j0}$. If $w_{\obfj}^{K,\gd_K}(\bfx)$ does not 
have singularity along $\caD_{j0}$ then clearly \eqref{iszero} 
is equal to 0. If it has singularity along $\caD_{j0}$ 
then \eqref{iszero} is the limit of
$\pm 2\pi iLi_{T_{\gd_K}^K(\bfj_{K-1})}
\Bigl(\bfx\big(T_{\gd_K}^K(\bfj_{K-1})\big)\Bigr)$ 
evaluated at $(\eps,\dots,\eps,0,\eps,\dots,\eps)$
where $0$ is at the $j$th place as  $\veps\to 0$ because $q(\veps)$
lies in $\cap_{i\ne j}\{x_i=\veps\}$ and does not enclose any 
irreducible component of $X_n$. This limit is 
equal to 0 by the induction assumption.

For brevity we drop the limit and 
$\veps$ in the rest of the proof. From Lemma~\ref{chen's}(iii)
\begin{multline*}
\left(\int_{qp}-\int_p\right)
\sum_{\substack{(\gd_1,\dots,\gd_K)\in \SS_K\\
\obfj\in\SS^K(m_1,\dots,m_n)}} \operatorname*{\bigsqcup}_{r=1}^K
w_{\obfj}^{r,\gd_r}\left(\bfx\Bigl(T_{\gd_{r+1}}^{r+1}
	\circ \cdots\circ T_{\gd_K}^K(\bfj_r) \Bigr)\right)\\
=\sum_{s=1}^{K}
\int_q\sum_{\substack{(\gd_1,\dots,\gd_K)\in \SS_K\\
\obfj\in\SS^K(m_1,\dots,m_n)}} \operatorname*{\bigsqcup}_{r=1}^s
w_{\obfj}^{r,\gd_r}\left(\bfx\Bigl(T_{\gd_{r+1}}^{r+1}
	\circ \cdots\circ T_{\gd_K}^K(\bfj_r) \Bigr)\right)\\
\cdot\int_p \operatorname*{\bigsqcup}_{r=s+1}^K
w_{\obfj}^{r,\gd_r}\left(\bfx\Bigl(T_{\gd_{r+1}}^{r+1}
	\circ \cdots\circ T_{\gd_K}^K(\bfj_r) \Bigr)\right).
\end{multline*}
We want to show that for each fixed $s$ the inner sum is zero. The
formula \eqref{Kcase} shows this is true for $s=K$. When $s=K-1$
we can divide the sum over products like
$$\int_q \phi_1\cdots \phi_{K-1}\int_p \phi_K$$
into sub-sums each one of which is produced by grouping all 
terms with $\phi_K=w_{\obfj}^{K,\gd_K}(\bfx)$
for some fixed $\gd_K$ and $\bfj_{K-1}$ which means that $\phi_K$ is
fixed. We see that the iterated integral of every sub-sum 
is $0$ by using equation~\eqref{ddd}
and applying \eqref{Kcase} with $K$ replaced
by $K-1$ and $\bfx$ replaced by 
$\bfx\big(T_{\gd_K}^K(\bfj_{K-1}) \big)$.

For any $1\le s\le n-2$ the sum over products like
$$\int_q \phi_1\cdots \phi_s\int_p \phi_{s+1}\cdots\phi_K$$
can be treated similarly by fixing $\phi_{s+1}\dots\phi_K$
first and then applying \eqref{Kcase} with $K=s$ and 
$\bfx$ replaced by $\bfx\big(T_{\gd_{s+1}}^{s+1}
	\circ \cdots\circ T_{\gd_K}^K(\bfj_s)\big)$.

This completes the proof of the theorem.
\end{proof}

\section{Multiple logarithms}
To study the mixed Hodge structure associated with
the multiple polylogarithms 
it is imperative that we resolve the monodromy of them.
In this section we carry this out for
multiple logarithm $\gemL_n(x_1,\dots,x_n)$.
We first provide a cleaner form of its analytic continuation.

Keeping the notation in the previous sections we have
$\SS(1,\dots,1)=\SS_n$ and $K=n$ for multiple logarithms. 
Though we can get the analytic continuation 
of the multiple logarithms by \eqref{e:genL} immediately, 
we actually have a cleaner expression in this special case. 

For any $\bfi=(i_1,\dots,i_n)\in \SS_n$ with $i_s=0$ we define
$$\text{pos}(\bfi,\bfi+\bfu_s)=s$$
as the position where the component is increased by 1.
For example $\text{pos}\big((1,0),(1,1)\big)=2$.
We define the position functions $f_n^1,\dots,f_n^n$ on
$\obfj\in \SS_n^n$ as follows:
$$f_n^1(\obfj)=1,\quad
f_n^t(\obfj)=\text{pos}\big(\bfj_{t-1},\bfj_t\big),
\text{ for \ } 2\le t\le n.$$ 
These functions tell us the places where the
increments occur in the queue of $\obfj$. Let
$$w_1(\bfx):=d\log\Bigl(\frac{1}{1-x_1}\Bigr);\quad
w_t(\bfx):=d\log\Bigl(\frac{1-x_{t-1}^{-1}}{1-x_t}\Bigr),
\text{ for \ }2\le t\le n.$$

\begin{prop} The multiple logarithm $\gemL_n(\bfx)$ is a 
multi-valued holomorphic function on $S_n$ and can be expressed by
\begin{equation} \label{formula}
\gemL_n(\bfx)=\sum_{\obfj=(\bfj_1,\dots,\bfj_n)\in\SS_n^n}
\int_{\zero}^{\bfx} w_{f_n^1(\obfj)}(\bfx(\bfj_1))
w_{f_n^2(\obfj)}(\bfx(\bfj_2))
\cdots w_{f_n^n(\obfj)}(\bfx(\bfj_n))
\end{equation}
\end{prop}
\begin{proof} Similar to the proof of Proposition \ref{genL}
and Theorem \ref{thm:mono}.
\end{proof}

We now turn to the monodromy of multiple logarithms.
\begin{lem} \label{s=n}
Let $p$ be a path from $\zero$ to $\bfx$ in $S_n$.
Let $q\in \pi_1(S_n,\bfx)$ be a loop turning around the component 
$\caD_{nn}=\{x_n-1=0\}$ only once but no other irreducible components 
of $X_n$
such that $\int_q d\log(1-x_n)=-2\pi i$. Then
\begin{equation*}
(\gTH(q)-\id) \gemL_n(\bfx)
=-2\pi i\gemL_{n-1}(x_1,\dots,x_{n-1}).
\end{equation*}
\end{lem}
\begin{proof} By Lemma \ref{pathindpt} we may assume that
$q$ is based at $\zero$ instead of $\bfx$.
We begin by moving the base of $p$ and $q$ to $\bfeps_n$ near $\zero$ 
and later we take the limit $\bfeps_n\to \zero$.

Let $p(\veps)$ and $q(\veps)$ be the corresponding
loop based at $\bfeps_n=(\veps,\cdots,\veps)\in D_n$.
By Lemma~\ref{indpath} we can
take the loop $q(\veps)$ in the two dimensional plane (over $\R$)
$x_1=\dots=x_{n-1}=\veps$ counterclockwise. By Lemma~\ref{chen's}(ii)
\begin{multline*}
(\int_{q(\veps)p(\veps)}-\int_{p(\veps)})\sum_{\obfj=
(\bfj_1,\dots,\bfj_n)
\in\SS_n^n}
w_{f_n^1(\obfj)}(\bfx(\bfj_1))
w_{f_n^2(\obfj)}(\bfx(\bfj_2))
\cdots w_{f_n^n(\obfj)}(\bfx(\bfj_n))\\
=\sum_{\obfj=(\bfj_1,\dots,\bfj_n)\in\SS_n^n}
\sum_{t=1}^{n} \int_{q(\veps)} w_{f_n^1(\obfj)}(\bfx(\bfj_1))
\cdots w_{f_n^t(\obfj)}(\bfx(\bfj_t))\\
\cdot \int_{p(\veps)}  w_{f_n^{t+1}(\obfj)}(\bfx(\bfj_{t+1}))
\cdots w_{f_n^n(\obfj)}(\bfx(\bfj_n)).
\end{multline*}

If $\bfj_1=\bfu_s$ for some $1\le s\le n-1$ then
$$w_1(\bfx(\bfj_1))=\frac{d(x_s\dots x_n)}{1-x_s\dots x_n}.$$
Thus for any $t=2,\cdots,n$
\begin{multline*}
\int_{q(\veps)} w_{f_n^1(\obfj)}(\bfx(\bfj_1))\cdots
w_{f_n^t(\obfj)}(\bfx(\bfj_t))\\
=\int_{q(\gep)}
-\log\left(\frac{1-\gep^{n-s}x_n}{1-\gep^{n-s+1}}\right)
w_{f_n^2(\obfj)}(\bfx(\bfj_2))\cdots
w_{f_n^t(\obfj)}(\bfx(\bfj_t))\to 0 \ \text{ as }\ \veps\to 0.
\end{multline*}

We now only need to look at those
$\obfj\in\SS_n^n$ with $\bfj_1=(0,\dots,0,1)$. Then
the last component of $\bfj_t$ is always 1 for $t=2,\dots, n$.
Thus the variable $\bfx(\bfj_t)=(\dots, x_n)$ has $t$ components
and $f_n^t(\bfj)\le t-1$ because $0$ cannot appear at the last
position of $\rho_{\bfj_t}(\bfj_{t-1})$.
Therefore,  $w_{f_n^t(\obfj)}(\bfx(\bfj_t))$ does not
involve the variable $x_n$ for $t=2,\dots, n$.

Suppose now $\phi_2,\dots,\phi_n$ are 1-forms on $S_n$ that
do not involve $x_n$. Then by Lemma~\ref{chen's}(ii)
$$\aligned
\ &\int_{q(\veps)p(\veps)}\frac{dx_n}{1-x_n} \phi_2\cdots 
\phi_n-\int_{p(\veps)}\frac{dx_n}{1-x_n} \phi_2\cdots  \phi_n\\
= &\int_{q(\veps)}\frac{dx_n}{1-x_n} \int_{p(\veps)} \phi_2\cdots \phi_n
+\sum_{t=2}^n\int_{q(\veps)}\frac{dx_n}{1-x_n} \phi_2\cdots \phi_t
\int_{p(\veps)}\phi_{t+1}\cdots \phi_n\\
=&-2\pi i \int_{p(\veps)} \phi_2\cdots \phi_n
+\sum_{t=2}^n(-1)^{t+1}\int_{q(\gep)^{-1}}
\phi_t\cdots \phi_2\frac{dx_n}{1-x_n}
\int_{p(\gep)} \phi_{t+1}\cdots \phi_n.
\endaligned$$
But none of $\phi_t$ for $t=2,\cdots,n$, involves variable $x_n$ by
assumption and therefore the inner most integral
$$\int_{(\veps,\cdots,\veps,\veps)} ^{(\veps,\cdots,\veps, x_n)} 
\phi_t=0$$
for any $x_n\in S_n$. Hence
$$\lim_{\veps\to 0}\Bigl(\int_{q(\veps)p(\veps)}-\int_{p(\veps)}\Bigr)
\frac{dx_n}{1-x_n} \phi_2\cdots \phi_n
=-2\pi i \int_p \phi_2\cdots \phi_n.$$
The lemma now follows from the one-to-one correspondence
$$\SS_{n-1}^{n-1}\longleftrightarrow \{\obfj\in \SS_n^n:
\bfj_1=\bfu_1\}.$$
\end{proof}

\begin{thm} \label{thm:monoDsn1}
Let $1\le s\le n$. Let $p$ be a path from $\zero$ to
$\bfx$ in $S_n$. Let $q_s\in \pi_1(S_n,\bfx)$ enclose
the component $\caD_{sn}=\{x_s\cdots x_n=1\}$ only once in $S_n$ but no
other irreducible components of $X_n$ such that 
$\int_{q_s} d\log(1-x_s\cdots x_n)=-2\pi i$. Then
\begin{equation*}
(\gTH(q_s)-\id) \gemL_n(\bfx)
=-2\pi i\gemL_{s-1}(x_1,\dots,x_{s-1})\cdot \gemL_{n-s}(\bfy(s))
\end{equation*}
where
$$\bfy(s)=\Bigl(\frac{1-x_sx_{s+1}}{1-x_s},\dots,
\frac{1-x_s\dots x_n}{1-x_s\dots x_{n-1}} \Bigr).$$
\end{thm}
\begin{proof} The case $s=n$ is proved by
Lemma~\ref{s=n}. We now prove the case $s=1$. The general
case will follow from these two cases by shuffle relations.

By Lemma \ref{indpath} we may  assume that the path $q$ lie
entirely in the (real) 2-dimensional plane $x_1=x_1$, $\dots$,
$x_{n-1}=x_{n-1}$. Suppose $s=1$ and $\bfy=\bfy(1)=(y_1,\dots,y_{n-1})$.
Let $\tilde{\SS}_n=\{\bfj=(j_1,\dots,j_n)\in \SS_n: j_1=1\}$ and
$$\tilde{\SS}_n^n=\{(\bfj_1,\dots, \bfj_n)\in \SS_n^n:
\bfj_t\in \tilde{\SS}_n \text{ for }t=1,\dots,n\}.$$
Then we can quickly find as in the proof of Lemma \ref{s=n} that
$$\aligned
(\int_{pq}-\int_p)  & \sum_{\obfj=(\bfj_1,\dots,\bfj_n)\in\SS_n^n}
w_{f_n^1(\obfj)}(\bfx(\bfj_1))
w_{f_n^2(\obfj)}(\bfx(\bfj_2))
\cdots w_{f_n^n(\obfj)}(\bfx(\bfj_n))\\
=&\int_q\sum_{\obfj=(\bfj_1,\dots,\bfj_n)\in\tilde{\SS}_n^n}
w_{f_n^1(\obfj)}(\bfx(\bfj_1))
w_{f_n^2(\obfj)}(\bfx(\bfj_2))
\cdots w_{f_n^n(\obfj)}(\bfx(\bfj_n))\\
=&(-1)^n\sum_{\obfj=(\bfj_1,\dots,\bfj_n)\in\tilde{\SS}_n^n}
\int_{q^{-1}} w_{f_n^n(\obfj)}(\bfx(\bfj_n))
\cdots w_{f_n^2(\obfj)}(\bfx(\bfj_2))
w_{f_n^1(\obfj)}(\bfx(\bfj_1))\\
=&-2\pi i\sum_{\obfj=(\bfj_1,\dots,\bfj_n)\in\tilde{\SS}_n^n}
\int_{(x_1,\dots,x_{n-1},(x_1\dots x_{n-1})^{-1})}^{(x_1,\dots,x_{n-
1},x_n)}
w_{f_n^2(\obfj)}(\bfx(\bfj_2))
\cdots w_{f_n^n(\obfj)}(\bfx(\bfj_n))
\endaligned$$
by Lemma \ref{chen's}(iii). Here the path in the integral
inside the last sum is a contractible path in $S_n$. Define
$$\aligned
\mu:\quad \SS_n &\longrightarrow \SS_{n-1} \\
(i_1,\dots,i_n) &\longmapsto (i_2,\dots,i_n)
\endaligned$$
and extend it to $\SS_n^n$ by mapping $(\bfi_1,\dots,\bfi_n)$
to $(\mu(\bfi_2),\dots,\mu(\bfi_n))$. This clearly induces a bijection
$\tilde{\SS}_n^n\leftrightarrow \SS_{n-1}^{n-1}$.
Now the case $s=1$ immediately follows 
from the following claim and the fact that
$\gemL_{n-1}(\bfy)$ evaluated at $(x_1,\dots,x_{n-1}, (x_1\dots
x_{n-1})^{-1})$ is 0. 
\begin{quote}
{\bf Claim.} For $2\le r\le s$
$$w_r\big(\bfx(\bfj_s)\big)=w_{r-1}\big(\bfy(\mu(\bfj_s))\big)$$.
\end{quote}
\noindent{\em Proof of the Claim.}
Let $\bfj_s=\sum_{i=1}^s \bfu_{t_i}$ where
$1=t_1<\cdots<t_s\le n$. Then
$\mu(\bfj_s)=\sum_{i=2}^s \bfu_{t_i-1}$ in $\SS_{n-1}$ and
$$\aligned
\bfx(\bfj_s)&=(x_1\cdots x_{t_2-1},x_{t_2}\cdots x_{t_3-1},
\dots,x_{t_s}\cdots x_n) \\
\bfy(\mu(\bfj_s))&=(y_{t_2-1}\cdots y_{t_3-2},y_{t_3-1}\cdots y_{t_4-2},
\dots, y_{t_s-1}\cdots y_{n-1}).
\endaligned$$
When $r=2$ we have
$$w_1\big(\bfy(\mu(\bfj_s))\big)=w_1(y_{t_2-1}\cdots y_{t_3-2})
=w_1\Bigl(\frac{1-x_1\cdots x_{t_3-1}}{1-x_1\cdots x_{t_2-1}} \Bigr)
=w_2\big(\bfx(\bfj_s)\big)$$
by the obvious identity
\begin{equation}\label{w2w1}
w_2(X,Y)= w_1\Bigl(\frac{1-XY}{1-X}\Bigr).
\end{equation}
By the same identity, when $r>2$ we have
$$w_r\big(\bfx(\bfj_s)\big)=w_2(x_{t_{r-1}}\cdots x_{t_r-1},
x_{t_r}\cdots x_{t_{r+1}-1})
=w_1\Bigl(\frac{1-x_{t_{r-1}}\cdots x_{t_{r+1}-1}}
{1-x_{t_{r-1}}\cdots x_{t_r-1}}\Bigr).$$
On the other hand
$$\aligned
w_{r-1}\big(\bfy(\mu(\bfj_s))\big)=&w_2(y_{t_{r-1}-1}\cdots y_{t_r-2},
y_{t_r-1}\cdots y_{t_{r+1}-2})\\
=&w_2\Bigl(\frac{1-x_1\cdots x_{t_r-1}}{1-x_1\cdots x_{t_{r-1}-1}},
\frac{1-x_1\cdots x_{t_{r+1}-1}}{1-x_1\cdots x_{t_r-1}}\Bigr)\\
\text{(by \eqref{w2w1} again)}\qquad =&w_1\left(
\Bigl(1-\frac{1-x_1\cdots x_{t_{r+1}-1}}{1-x_1\cdots x_{t_{r-1}-1}}
\Bigr)/
\Bigl(1-\frac{1-x_1\cdots x_{t_r-1}}{1-x_1\cdots x_{t_{r-1}-1}}\right)\\
=&w_1\Bigl(\frac{x_1\cdots x_{t_{r+1}-1}-x_1\cdots x_{t_{r-1}-1}}
{x_1\cdots x_{t_r-1}-x_1\cdots x_{t_{r-1}-1}}\Bigr)\\
=&w_r\big(\bfx(\bfj_s)\big).
\endaligned$$
The claim now is proved.

\bigskip

Let $2\le s\le n-1$ and $\bfy=\bfy(s)$. By Lemma~\ref{indpath} we can
move our base point of the loop $q_s$ from $\bfx$ to $\bveps_n$ close 
to $\zero$. Let $q(\veps)$ be a loop around $\caD_{sn}$ based at $\bveps_n$ 
lying entirely in $x_1=\cdots x_{n-1}=\veps$ such that 
$\int_{q(\veps)} d\log(1-x_s\cdots x_n)=-2\pi i$. 
Let $p(\veps)$ be a path from $\bveps_n$ to $\bfx$ in $S_n$.
Then as in Lemma~\ref{s=n} we find that
\begin{align}
\lim_{\veps\to 0}(\int_{q(\veps)p(\veps)}-\int_{p(\veps)}) & 
\sum_{\obfj=(\bfj_1,\dots,\bfj_n)\in\SS_n^n}
w_{f_n^1(\obfj)}(\bfx(\bfj_1))
w_{f_n^2(\obfj)}(\bfx(\bfj_2))
\cdots w_{f_n^n(\obfj)}(\bfx(\bfj_n))\notag \\
=& -2\pi i\sum_{\obfj\in\SS_n^n,\ \bfj_1=\bfu_s}
\lim_{\veps\to 0}\int_{p(\veps)} w_{f_n^2(\obfj)}(\bfx(\bfj_2))
\cdots w_{f_n^n(\obfj)}(\bfx(\bfj_n)) \label{sgen}
\end{align}
For any $\bfj\in \SS_n^n$ with $\bfj_1=\bfu_s$ we define the type
of $\bfj_t$ ($t\ge 2$) by
$$\text{Type}(\bfj_t)=\begin{cases}
\text{(I)} \quad &\text{ if }\bfj_t=\bfj_{t-1}+\bfu_m \text{ and }m<s; 
\\
\text{(II)} \quad &\text{ if }\bfj_t=\bfj_{t-1}+\bfu_m\text{ and } m>s.
\end{cases}$$
Note that $m=f_n^t(\obfj)$. Now each term in the sum 
\eqref{sgen} has the form $\phi_2\cdots \phi_n$ with $s-1$ of
them being of type (I) and $n-s$ of them being of type (II).
If $\phi_t=w_{f_n^t(\obfj)}(\bfx(\bfj_t))$ is
of type (I) then $f_n^t(\obfj)<s$ and the $s$-th component of
$\bfj_t$ is 1 which implies that $\phi_t$ does not involve the variables
$x_s,\dots,x_n$. Similar argument applies to type (II) 1-forms which
does not involve the variables $x_1,\dots, x_{s-1}$. Furthermore, 
altogether there are $(n-1)!$ terms in the sum of \eqref{sgen}
which are all different from
each other. An easy computation using Lemma \ref{s=n} for case $s=n$
and the above claim in the proof of case $s=1$
then yields that each term in the sum is a term in the expansion of
the right hand side of (setting $\bfx'=(x_1,\dots,x_{s-1})$)
\begin{equation}\label{prodd}
\aligned
\gemL_{s-1}(\bfx')\gemL_{n-s}(\bfy)
=&\lim_{\veps\to 0}\int_{p(\veps)}\sum_{\obfi\in\SS_{s-1}^{s-1}}
w_{f_{s-1}^1(\obfi)}(\bfx'(\bfi_1))
\cdots w_{f_{s-1}^{s-1}(\obfi)}(\bfx'(\bfi_{s-1})) \\
& \cdot \lim_{\veps\to 0}\int_{p(\veps)}\sum_{\obfk\in\SS_{n-s}^{n-s}}
w_{f_{n-s}^1(\obfk)}(\bfy(\bfk_1))
\cdots w_{f_{n-s}^{n-s}(\obfk)}(\bfy(\bfk_{n-s}))
\endaligned
\end{equation}
after the shuffle relation of Lemma \ref{chen's}(iv) is applied to it.
We have used the fact that in the computation of case $s=1$, we can
move the base from $\bfx$ to $\bveps_n$ and then take $\veps\to 0$.
This provides the expression of the function $\gemL_{n-s}(\bfy)$ in the
above.

Now on the right hand side of \eqref{prodd} there are $(s-1)!(n-s)!$ 
terms
each of which produces exactly ${n-1\choose s-1}$ terms by
the shuffle relation. Therefore the right hand side of \eqref{prodd}
will produce exactly $(n-1)!$ terms after the
shuffle relation is applied.

This completes the proof of the theorem.
\end{proof}

The next result shows that the monodromy produced by Theorem
\ref{thm:monoDsn1} is the only kind of monodromy of the multiple
logarithms.

\begin{prop}\label{prop:monoDij0}
The monodromy of $\gemL_n(\bfx)$ about 
$\caD_{ii}=\{x_i=1\},$ $1\le i<n$, and
$\caD_{ij}=\{x_i\cdots x_j=1\}$, $1\le i<j<n,$ is trivial.
\end{prop}
\begin{proof} Let $p\in \pi_1(S_n,\bfx)$ be a path from $\zero$ to
$\bfx$ and $q\in \pi_1(S_n,\zero)$ be a loop which encloses 
$\caD_{ij}$ only once but no other irreducible component of $X_n$.

It suffices to look at the 1-forms in \eqref{formula} which
have singularities along $\caD_{ij}$, $1\le i\le j<n$. Suppose
$\obfj=(\bfj_1,\dots,\bfj_n)\in \SS_n^n$ and $\bfj_s=(j_1,\dots,j_n)$.
Because $\bfj_s$ has depth $s$ we may assume that
$j_{t_1},\dots,j_{t_s}$ are the only nonzero components of $\bfj_s$.
Let $\bfy=\bfx(\bfj_s)=(y_1,\dots,y_s)$ then
$$w_r(\bfy)=\frac{dy_r}{1-y_r}-\frac{dy_{r-1}}{y_{r-1}}
	-\frac{dy_{r-1}}{1-y_{r-1}}, \quad 1\le r\le s,$$
where the last two terms do not appear if $r=1$. This 1-form
has singularity along $\caD_{ij}$, $1\le i\le j<n$, if $t_r=i$ and
$t_{r+1}=j+1$. If this is the case then 
$$w_{r+1}(\bfy)=\frac{dy_{r+1}}{1-y_{r+1}}
-\frac{dy_r}{y_r}-\frac{dy_r}{1-y_r}$$ 
has singularity along $\caD_{ij}$ too. These two 1-forms correspond to
the following choices of $\bfj_{s-1}$ respectively:
$${\setlength\arraycolsep{1pt}
\begin{array}{rllll}
\bfj_{s-1}'=\bfj_s-\bfu_i=
(\dots,0,&1,0,\dots,0,&0,0,\dots,0,&1,0,\dots,0,&1,0,\dots).\\
\     &\uparrow  &\uparrow     &\uparrow &\uparrow  \\
\ &\hskip-7pt  t_r\text{-th}& \hskip-5pti\text{-th} & \hskip-14pt
(j+1)\text{-st}& \hskip-10pt t_{r+1}\text{-th} \\
\     &\downarrow  &\downarrow  &\downarrow  &\downarrow  \\
\bfj_{s-1}''=\bfj_s-\bfu_{j+1}=
(\dots,0,&1,0,\dots,0,&1,0,\dots,0,&0,0,\dots,0,&1,0,\dots).
\end{array}}$$
Note $j<n$ by assumption  so $\bfu_{j+1}$ makes sense. Thus 
for each fixed $\bfj_s$ with $t_r=i$ and $t_{r+1}=j+1$ we can
regroup all terms in \eqref{formula} which have singularity
along $\caD_{ij}$ into sub-sums of only two terms, one for
$$\obfj'=(\bfj_1,\dots,\bfj_{s-2},\bfj_{s-1}',\bfj_s,\bfj_
{s+1},\dots,\bfj_n)$$
and the other for
$$\obfj''=(\bfj_1,\dots,\bfj_{s-2},\bfj_{s-1}'',\bfj_s,\bfj_
{s+1},\dots,\bfj_n)
$$
with arbitrarily fixed
$\bfj_1,\dots,\bfj_{s-2},\bfj_{s+1},\dots,\bfj_n$, where  
we then must have $\bfj_{s-2}=\bfj_s-\bfu_i-\bfu_{j+1}$.
For such two terms we have by Lemma  \ref{chen's}(ii)
\begin{multline*}
(\gTH(q)-\id)\int_p\sum_{\obfj=\obfj', \obfj''}
 \operatorname*{\bigsqcup}_{r=1}^{n} w_{f_n^r(\obfj)}\big(\bfx(\bfj_r)
\big)
=\int_p \operatorname*{\bigsqcup}_{r=s+1}^{n} 
w_{f_n^r(\obfj)}\big(\bfx(\bfj_r)\big)\\
\cdot \int_q \operatorname*{\bigsqcup}_{r=1}^{s-2}
	w_{f_n^r(\obfj)}\big(\bfx(\bfj_r)\big)
	\Bigl[w_{f_n^{s-1}(\obfj')}\big(\bfx(\bfj_{s-1}')\big)
	w_r\big(\bfx(\bfj_s)\big)
	+w_{f_n^{s-1}(\obfj'')}\big(\bfx(\bfj_{s-1}'')\big)
	w_{r+1}\big(\bfx(\bfj_s)\big)\Bigr].
\end{multline*}
Here we have taken the liberty to drop the primes when it is the same to
write $'$ or $''$. We see that the second iterated integral is zero 
because
$$\aligned
w_{f_n^{s-1}(\obfj')}(\bfx(\bfj_{s-1}'))=&
w_r(\dots,x_{t_{r-1}}\cdots x_j,x_{j+1}\cdots x_{t_{r+1}},\dots)\\
w_{f_n^{s-1}(\obfj'')}(\bfx(\bfj_{s-1}''))=&
w_r(\dots,x_{t_{r-1}}\cdots x_{i-1},x_i\cdots x_{t_{r+1}},\dots)
\endaligned$$
coincide along $\caD_{ij}$. This finishes the proof of
the proposition.
\end{proof}

We end our paper with a result which will be used in the computation
of the mixed Hodge structures associated with multiple logarithms.
\begin{prop} \label{monoD}
Let $n>1$. For any $1\le a<b\le n$ set $F_{aa}=1$ and
$$F_{ab}(\bfx)=\gemL_{b-a} \left(\frac{1-x_ax_{a+1}}{1-x_a},\cdots,
\frac{1-x_a\cdots x_b}{1-x_a\cdots x_{b-1}}\right).$$
Let $1\le j< n$ and $q_{j0}\in \pi_1(S_n,\bfx)$
(resp. $1\le j< n$ and $q_{1j}$,  $2\le j\le n$ and $q_{jn}$) 
be a loop turning around the component $\caD_{j0}=\{x_j=0\}$ 
(resp. $\caD_{1j}=\{x_1\cdots x_j=1\}$, 
resp. $\caD_{jn}=\{x_j\cdots x_n=1\}$),
only once but no other irreducible components of $X_n$
such that $\int_{q_{j0}}dx_j/x_j=2\pi i$ 
(resp. $\int_{q_{1j}} d\log(1-x_1\cdots x_j)=2\pi i$,
resp. $\int_{q_{jn}}d\log(1-x_j\cdots x_n)=2\pi i$). Then
$$\aligned
(\gTH(q_{j0})-\id)F_{1n}(\bfx)=&
	-2\pi i \sum_{s=j}^{n-1} F_{1s}(\bfx)F_{s+1,n}(\bfx),\\
(\gTH(q_{1j})-\id)F_{1n}(\bfx)=&
	\phantom{-} 2\pi i F_{1,j}(\bfx)F_{j+1,n}(\bfx),\\	
(\gTH(q_{jn})-\id)F_{1n}(\bfx)=&-2\pi i F_{1,{j-1}}(\bfx)F_{jn}(\bfx), 
\endaligned$$
where $\gTH(q)$ denotes the action of $q\in \pi(S_n,\bfx)$.
\end{prop}
\begin{proof} The proposition follows from the monodromy property of
$\gemL_n(\bfx)$. We only prove the result for $\caD_{j0}$ because
the proof is exactly the same for $\caD_{1j}$ and $\caD_{jn}$.

By Theorem \ref{thm:monoDsn1} we know that the monodromy of 
$\gemL_n(\bfx)$ around $x_s\cdots x_n=1$ is given by 
$-2\pi i\gemL_{s-1}(x_1,\dots,x_{s-1}) F_{sn}(\bfx)$. Let
$y_s=\frac{1-x_1\cdots x_{s+1}}{1-x_1\cdots x_s}$ and
$\bfy=(y_1,\dots,y_{n-1})$. On $\caD_{j0}$ all of 
$y_s\cdots y_{n-1}=\frac{1-x_1\cdots x_n}{1-x_1\cdots x_s}$, 
$j\le s<n$, are equal to 1. Therefore
the monodromy of $F_{1n}(\bfx)=\gemL_{n-1}(\bfy)$ 
about $\caD_{j0}$ is
the sum 
$$-2\pi i\sum_{s=j}^{n-1} \gemL_{s-1}(y_1,\dots,y_{s-1}) F_{s,n-1}(\bfy)
=-2\pi i \sum_{s=j}^{n-1} F_{1s}(\bfx)F_{s+1,n}(\bfx).$$
This concludes our proposition and the paper.
\end{proof}

\bigskip

\noindent
{\em Address:} Department of Mathematics, University of Pennsylvania, 
PA 19104, USA

\noindent
{\em Email:} jqz@math.upenn.edu

\end{document}